\documentclass[12pt,reqno]{amsart}  
\usepackage{enumerate}
\usepackage{url} 

\usepackage{ifthen} 

\usepackage[italian,english]{babel}
\usepackage{comment} 
\usepackage{fancyhdr}

\renewcommand{\qedsymbol}{$\square$}

\newcommand{\dx}{\mathrm{d}}
\newcommand{\eps}{\varepsilon}
\newcommand{\R}{\mathbb{R}}

\def \fine {{\hfill \qedsymbol}}

\newtheoremstyle{slantedtheorems}
{10pt}
{6pt}
{\slshape}
{}
{\bfseries}
{.}
{.5em}
{\thmname{#1} \thmnumber{#2}\thmnote{ (#3)}}
 
 \theoremstyle{slantedtheorems}

  \newtheoremstyle{Boldtheorems}
  {10pt}
  {6pt}
  {\bfseries}
  {}
  {\bfseries}
  {.}
  {.5em}
  {\thmname{#1} \thmnumber{#2}\thmnote{ (#3)}}
  
\theoremstyle{Boldtheorems}

\allowdisplaybreaks

\title[]{An extension of the pair-correlation conjecture \\ and applications}
\author[]{A.~Languasco, A.~Perelli \lowercase{and} A.~Zaccagnini}

\numberwithin{equation}{section}

\begin{document}

\maketitle

\vskip.5cm
{\bf Abstract.} We study an extension of Montgomery's pair-correlation conjecture and its relevance in some problems on the distribution of prime numbers.



\section{Introduction}

\smallskip
Throughout the paper we assume the Riemann Hypothesis (RH) for the Riemann zeta function $\zeta(s)$. Let
\[
F(X,T) = 4 \sum_{-T\leq\gamma,\gamma'\leq T}
\frac{X^{i(\gamma-\gamma')}}{4+(\gamma-\gamma')^2},
\] 
where $\gamma,\gamma'$ run over the imaginary part of the non-trivial zeros of $\zeta(s)$, be the pair-correlation function of the zeta zeros.  In his famous paper \cite{Mon/1973}, Montgomery conjectured that
\begin{equation}
\label{1-1}
F(X,T) \sim \frac{T}{\pi}\log T
\end{equation}
as $X\to\infty$, uniformly for $X^\eps \leq T \leq X$ for any fixed $0<\eps <1$. Moreover, the essentially equivalent assertion in \cite{Mon/1973}, that for any fixed $0<\alpha<\beta$ and $T\to\infty$
\[
\frac{\pi}{T\log T} \sum_{\substack{-T\leq\gamma,\gamma'\leq T\\ \frac{2\pi\alpha}{\log T} \leq \gamma-\gamma'\leq \frac{2\pi\beta}{\log T}}} 1 \sim \int_\alpha^\beta \Big(1-\big(\frac{\sin\pi u}{\pi u}\big)^2\Big) \dx u,
\]
raised an interesting connection with the Random Matrix Theory; see the excellent survey by Conrey \cite{Con/2001}.  A sharp numerical confirmation of the latter conjecture was given by Odlyzko \cite{Odl/1987},\cite{Odl/1989}; see also \cite{Odl/web}. 

\medskip
The pair-correlation conjecture \eqref{1-1}, as well as its variants, found several applications in prime number theory. Among the vast literature on this subject we recall here the papers by Gallagher $\&$ Mueller \cite{Ga-Mu/1978}, Heath-Brown \cite{HB/1982}, Goldston $\&$ Montgomery \cite{Go-Mo/1987} and Rudnick $\&$ Sarnak \cite{Ru-Sa/1996}. In particular, writing
\[
J(X,h) = \int_X^{2X} (\psi(x+h)-\psi(x)-h)^2 \dx x \quad \text{and} \quad \psi(x) = x + \Delta(x),
\]
Heath-Brown \cite{HB/1982} showed that assuming RH and various forms of the pair-correlation conjecture, weaker than \eqref{1-1} concerning the size and/or the uniformity, one could get sharp bounds for the error term $\Delta(x)$, for the mean-square of primes in short intervals $J(X,h)$ and for the difference between consecutive primes $p_{n+1}-p_n$. The ideas and techniques in \cite{HB/1982} were further developed by several authors. We mention here the papers by Heath-Brown $\&$ Goldston \cite{HB-Go/1984} and Liu $\&$ Ye \cite{Li-Ye/2001}, the latter assuming RH and a version of the pair-correlation conjecture given by Rudnick $\&$ Sarnak \cite{Ru-Sa/1996}. We also mention the earlier paper by Mueller \cite{Mue/1981} on $p_{n+1}-p_n$. Moreover, we recall that Gonek \cite{Gon/1993} proposed a conjecture related to Landau's explicit formula and showed that such a conjecture has very strong consequences on the distribution of primes in short intervals. We shall enter the contents of the above quoted papers later on, while comparing to our results.

\medskip
It turns out that the applications of the pair-correlation conjecture contained in some of the above papers, along with further results, can be framed in a unified way as consequences of suitable assumptions on the extended pair-correlation function
\[
F(X,T,\tau) = 4 \sum_{-T\leq\gamma,\gamma'\leq T} \frac{X^{i(\gamma-\gamma')}}{4+\tau^{2}(\gamma-\gamma')^2},
\]
where $0\leq \tau \leq 1$. Note that at the extreme $\tau$-values we have $F(X,T,1) = F(X,T)$ and $F(X,T,0) = |\Sigma(X,T)|^2$, where
\begin{equation}
\label{1-2}
\Sigma(X,T) = \sum_{|\gamma|\leq T} X^{i\gamma}
\end{equation}
is the exponential sum over zeros appearing in Landau's explicit formula, see Gonek \cite{Gon/1993}. We first give an interpretation of $F(X,T,\tau)$ as a suitable pair-correlation function. Indeed, it is clear that the numbers $\widetilde{\rho}=\tau/2 + i\tau\gamma$ are the non-trivial zeros of the function $Z(s) = \zeta(s/\tau)$, and
\begin{equation}
\label{1-3}
F(X,T,\tau) = 4 \sum_{-\tau T\leq \tau\gamma,\tau\gamma'\leq \tau T} \frac{(X^{1/\tau})^{i(\tau\gamma-\tau\gamma')}}{4+(\tau\gamma-\tau\gamma')^2}.
\end{equation}
Moreover, given a rather general $L$-function $L(s)$ (namely, belonging to the Selberg class of $L$-functions, see e.g. Kaczorowski $\&$ Perelli \cite{Ka-Pe/1999b}) and assuming the Riemann Hypothesis for $L(s)$, Murty $\&$ Perelli \cite{Mu-Pe/1999} investigated the pair-correlation function $F_L(X,T)$ which, with a slightly different normalization with respect to \cite{Mu-Pe/1999}, is defined in analogy with $F(X,T)$ as
\begin{equation}
\label{1-4}
F_L(X,T) = 4 \sum_{-T\leq\widetilde{\gamma},\widetilde{\gamma}'\leq T}
\frac{(X^d)^{i(\widetilde{\gamma}-\widetilde{\gamma}')}}{4+(\widetilde{\gamma}-\widetilde{\gamma}')^2}.
\end{equation}
Here the non-trivial zeros of $L(s)$ are of the form $\widetilde{\rho} = 1/2 + i\widetilde{\gamma}$, and $d$ is the degree of $L(s)$ (i.e. twice the sum of the coefficients $\lambda_j$ in the $\Gamma$-factors $\Gamma(\lambda_js+\mu_j)$ appearing in its functional equation; see \cite{Ka-Pe/1999b}). Since $\zeta(s)$ has degree 1, we have $F_\zeta(X,T)=F(X,T)$. The functional equation of $Z(s)$, obtained at once from the one of $\zeta(s)$, suggests that we may heuristically regard $Z(s)$ as an $L$-function of degree $d$ and conductor $q$ given by
\begin{equation}
\label{1-5}
d=\frac{1}{\tau} \hskip1.5cm q=\big(\frac{1}{\tau}\big)^{\frac{1}{\tau}}
\end{equation}
(see again \cite{Ka-Pe/1999b} for the definition of conductor), although $Z(s)$ does not belong to the Selberg class if $\tau\neq 1$ (e.g. since it has a pole at $s=\tau$). Hence in view of \eqref{1-3}, \eqref{1-4} and \eqref{1-5} we have
\begin{equation}
\label{1-6}
F(X,T,\tau) = F_Z(X,\tau T).
\end{equation}
The pair-correlation conjecture in Murty $\&$ Perelli \cite{Mu-Pe/1999}, namely
\[
F_L(X,T) \sim \frac{dT}{\pi} \log T
\]
as $X\to\infty$, uniformly for $X^\eps \leq T \leq X$ for any fixed $0<\eps<1$, is meant for a given $L$-function, hence in particular for fixed degree $d$. In view of the Riemann-von Mangoldt formula for the the Selberg class, a reasonable $(d,q)$-uniform version of this conjecture is 
\begin{equation}
\label{1-7}
F_L(X,T) \sim \frac{dT}{\pi} \log q^{1/d}T
\end{equation}
in the same range of uniformity for $X$ and $T$. Hence, in view of \eqref{1-5}, \eqref{1-6} and \eqref{1-7} we may expect that
\begin{equation}
\label{1-8}
F(X,T,\tau) \ll T\log T
\end{equation}
uniformly for $X^\eps \leq \tau T \leq X$, and possibly even in wider ranges for $T$ and $\tau$.  We refer to \eqref{1-7} as the general Montgomery conjecture.

\medskip
We recall that the asymptotic behavior of $F(X,T)$ in \eqref{1-1} is expected to change around $T=X$ and that, under RH, Montgomery \cite{Mon/1973} (and Goldston \cite{Gol/1981}) detected the behavior of $F(X,T)$ when $T\geq X$. An analog of the latter result is proved in Murty $\&$ Perelli \cite{Mu-Pe/1999} for the functions $F_L(X,T)$; in such case the assumption $T\geq X^d$ is required. However, in our case the function $Z(s)$ is a $\tau$-rescaling of $\zeta(s)$, and for this reason we may expect that the role of $X$ is played essentially by $X^\tau$. In the companion paper \cite{L-P-Z/tauPC-II} we obtain, assuming certain hypotheses, the asymptotic behavior of $F(X,T,\tau)$ in suitable ranges of $T$ and $\tau$, when $T\geq X$. Note that the trivial bound for $F(X,T)$ is $O(T\log^2T)$ uniformly for $X,T\geq 2$, hence only $\log T$ worse than the expected size (when $X^\eps\leq T\leq X$), while the trivial bound for $F(X,T,\tau)$ is for small $\tau$ much worse than the expected size (again if $T\leq X$). Indeed, we have
\begin{equation}
\label{1-9}
F(X,T,\tau) \ll \min \Bigl(T; \frac{1}{\tau}\Bigr)\ T \log^2 T  
\end{equation}
uniformly for $X,T\geq 2$ and $0\leq \tau\leq 1$. This is obvious for $\tau\leq 1/T$, while in the opposite case we have
\[
F(X,T,\tau) \ll T\log^2T \sum_{m=0}^{T}  \frac{1}{1+\tau^2m^2} \ll \frac{1}{\tau} T\log^2T.
\]
In view of \eqref{1-8} and \eqref{1-9}, in this paper we also consider the slightly more conservative pair-correlation conjecture
\begin{equation}
\label{1-10}
F(X,T,\tau) \ll TX^\eps
\end{equation}
for every $\eps>0$, with $T$ and $\tau$ in suitable ranges.

\medskip
The extended pair-correlation function $F(X,T,\tau)$ already appears in Heath-Brown $\&$ Goldston \cite{HB-Go/1984} (as $G_\beta(X,T)$, $\beta=1/\tau$, in their notation) and is used there to study the size of $p_{n+1}-p_n$ under Montgomery's conjecture \eqref{1-1} with $T$ in a certain range. To this end, $F(X,T,\tau)$ is related to $F(X,T)$ (see Lemma 2 in \cite{HB-Go/1984}), and the required bounds on $F(X,T,\tau)$ are obtained from the asymptotic formula for $F(X,T)$. However, the information obtained on $F(X,T,\tau)$ in that way reflects its expected order of magnitude only in a rather short range of $\tau$ close to 1. Clearly, that approach allows to obtain bounds for $F(X,T,\tau)$ from asymptotic formulae with remainder for $F(X,T)$, say of the form considered in our recent paper \cite{L-P-Z/2012} (roughly, an error of type $T\eps(T)$ in input gives an error of type $T\eps(T)/\tau^2$ in output). We remark that, in \cite{HB-Go/1984} and in subsequent works on the subject, the function $F(X,T,\tau)$ is used essentially as an intermediate tool in order to get information on the distribution of primes from suitable forms of Montgomery's pair-correlation conjecture. Our aim in this paper is, instead, to show how $F(X,T,\tau)$ directly controls the distribution of primes, especially in short intervals, depending on the range of uniformity in $T$ and $\tau$ allowed  in \eqref{1-10} or \eqref{1-8}. Actually, we mainly focus on the distribution of primes in short intervals of length $h\geq X^\eps$, with $\eps>0$ arbitrarily small, under the following hypotheses on the uniformity ranges for $T$ and $\tau$ where \eqref{1-10} is assumed to hold.

\medskip
{\bf Hypothesis $H(\eta)$.} {\sl Assume RH, let $0\leq \eta<1$ be fixed and let $X\to\infty$. If $0<\eta<1$ then \eqref{1-10} holds uniformly for $X^\eta \leq T \leq X$ \and $X^\eta/T \leq \tau \leq 1$; if $\eta=0$ then for every small $\eps>0$, \eqref{1-10} holds uniformly for $X^\eps \leq T \leq X$ \and $0 \leq \tau \leq 1$.}

\medskip
{\bf Remark 1.} First we note that for any fixed $0< \eta<1$, the uniformity range of $\tau T$ in hypothesis $H(\eta)$ is contained in $[X^\eps, X]$, which in view of \eqref{1-6}-\eqref{1-8} may be regarded as a plausible range of a $d$-uniform version of the general Montgomery conjecture. When $\eta=0$, $\tau T$ may become $\leq X^\eps$. However, as we remarked earlier, Gonek \cite{Gon/1993} made a conjecture on the size of the exponential sum $\Sigma(X,T)$ in \eqref{1-2}, namely
\begin{equation}
\label{1-11}
\Sigma(X,T) \ll T X^{-1/2+\eps} + T^{1/2}X^\eps
\end{equation}
for $X,T\geq2$. When $T\leq X$ the second term dominates, and hence in view of \eqref{1-11} and $F(X,T,0)=|\Sigma(X,T)|^2$ we may expect that for $X^\eps\leq T \leq X$ and $0\leq \tau\leq 1$
\[
F(X,T,\tau) \ll |\Sigma(X,T)|^2 \ll TX^\eps.
\]
Actually, by Lemma 5 below we have that for $\tau T\leq X^\eps$
\[
F(X,T,\tau) \ll T^{1+\eps} + X^\eps \max_{t\leq T}|\Sigma(X,t)|^2,
\]
thus \eqref{1-10} follows from Gonek's conjecture in this range. Summarizing, we may say that our hypothesis on the size of $F(X,T,\tau)$ is supported by the general Montgomery conjecture in the range $X^\eps\leq \tau T\leq X$, and by Gonek's conjecture in the remaining range $\tau T\leq X^\eps$. Actually, $H(0)$ is equivalent to Gonek's conjecture restricted to the range $X^\eps\leq T\leq X$. \fine

\medskip
Assuming RH and the pair-correlation conjecture in the form $F(X,T)\ll T\log T$ uniformly for $X^\eps\leq T \leq X$, Heath-Brown \cite{HB/1982} implicitly proved (see the proof of Theorem 2 in \cite{HB/1982}) that $J(X,h) \ll hX\log(2X/h)$ uniformly for $1 \leq h \leq X^{1-\eps}$. Later, Goldston $\&$ Montgomery \cite{Go-Mo/1987} proved that the conjectural asymptotic behavior of $F(X,T)$ implies that of $J(X,h)$, and vice versa. Assuming only RH one has $J(X,h) \ll hX\log^2(2X/h)$ uniformly for $1\leq h \leq X$, see Selberg \cite{Sel/1942} and Saffari $\&$ Vaughan \cite{Sa-Va/1977}. Let
\[
J(X,Y,h) = \int_X^{X+Y} \bigl( \psi(x + h) - \psi(x) - h \bigr)^2  \dx x.
\]
We expect that, as $X\to\infty$,
\begin{equation}
\label{1-12}
J(X,Y,h) \ll hYX^\eps
 \end{equation}
uniformly in suitable ranges of $h$ and $Y$, for every $\eps>0$. We start with a result of Heath-Brown's type for $J(X,Y,h)$, under the assumption of hypothesis $H(\eta)$; we deal in \cite{L-P-Z/tauPC-II} with the analog of Goldston $\&$ Montgomery \cite{Go-Mo/1987}, i.e. with the equivalence between the asymptotic behavior of $J(X,Y,h)$ and $F(X,T,\tau)$.

\medskip
{\bf Theorem 1.} {\sl Assume $H(\eta)$. If $0<\eta<1$ then \eqref{1-12} holds uniformly for $1 \leq h \leq X^{1-\eta}$ and $hX^\eta\leq Y\leq X$. If $\eta=0$ then \eqref{1-12} holds uniformly for $1 \leq h \leq X^{1-\eps}$ and $0\leq Y \leq X$, for every small $\eps>0$.}

\medskip
By similar arguments one can get variants of Theorem 1, for example where both the hypothesis and the result are localized, or replacing $X^\eps$ (both in input and in output) by suitable powers of $\log X$. Moreover,  a standard consequence of Theorem 1 is that 
\[
\psi(x+h)-\psi(x)=h+O(h^{1/2}x^\eps)
\] 
for almost all $x\in[X,X+Y]$, with $h$ and $Y$ as in Theorem 1. 

\medskip
Hypothesis $H(\eta)$ and Theorem 1 have further consequences on the distribution of primes in short intervals. For example, using in addition the inertia property of the function $\psi(x)$, from Theorem 1 we can deduce results on $\psi(x+h)-\psi(x)$ valid for all large $x$. On the other hand, hypothesis $H(\eta)$ can be used in a direct way for the same purpose. It is interesting to note that both approaches have advantages over the other, and the output of the two approaches is given by

\medskip
{\bf Theorem 2.} {\sl Assume $H(\eta)$ and let $\eps>0$ be arbitrarily small. If $0<\eta<1/2-5\eps$ then
\[ 
\psi(x+h)-\psi(x) = h +
\begin{cases}
O(h^{2/3}x^{\eta/3+\eps}) &  \text{for} \  x^{\eta+5\eps}\leq h\leq x^{1/2} \\
O(h^{1/3} x^{1/6+\eta/3+\eps}) & \text{for} \  x^{1/2}\leq h\leq x^{1-\eta}.
\end{cases} 
\]
If $\eta=0$ then $\psi(x+h)-\psi(x)= h+O( h^{1/2} x^\eps)$ for $x^{3\eps} \leq h \leq x^{1 - \eps}$.} 

\medskip
{\bf Remark 2.} Let $0<\eta<1/2-5\eps$. The first estimate is obtained from Theorem 1 by means of the inertia property, while the second follows by a direct application of $H(\eta)$. Both estimates are non-trivial in the stated ranges for $h$, in the sense that the error terms are $O(\min(h^{1-\eps},x^{1/2+\eps}))$ and the ranges are non-empty (recall that the second bound holds under RH). However, we believe that the direct approach should always give a better result, namely
\[
\psi(x+h)-\psi(x)=h+O(h^{1/2}x^{\eta/2+\eps})
\]
for $x^{\eta+4\eps}\leq h\leq x^{1-\eta}$ and $0<\eta<1/2-2\eps$. At present we cannot get such a bound because of the maximum over $t$ in Lemma 6 below, which forces the use of trivial bounds in a certain range, thus giving the weaker error term in Theorem 2. Note that there is a discontinuity between the results in the cases $\eta>0$ and $\eta=0$, due to an analogous discontinuity in hypothesis $H(\eta)$. Note also that hypothesis $H(\eta)$ with $\eta>0$ has two distinct consequences on primes in short intervals, one in mean-square (Theorem 1) and the other for each large $x$ (Theorem 2). The situation changes when $\eta=0$, since it is easily seen that in this case Theorem 1 is equivalent to Theorem 2. This is in agreement with the observation about Theorem 4 of Gonek \cite{Gon/1993} at the end of Remark 1, and indeed our bound for $\psi(x+h)-\psi(x)-h$ under $H(0)$ is the same as Gonek's bound under his conjecture. Finally, in the direction of Mueller \cite{Mue/1981}, Heath-Brown \cite{HB/1982}, Heath-Brown $\&$ Goldston \cite{HB-Go/1984} and Liu $\&$ Ye \cite{Li-Ye/2001}, bounds for $p_{n+1}-p_n$ of type $p_n^{1/2}f(p_n)$, with suitable $f(x)$'s, can be obtained from the estimate $F(X,T,\tau)\ll T\log T$ with suitable uniformity ranges of $T$ and $\tau$. In this case $T$ will be around $X^{1/2}$ and $\tau$ around a negative power of $\log X$. \fine

\bigskip
Now we turn to the problem of estimating $\Delta(x) = \psi(x)-x$ under suitable pair-correlation hypotheses. The bound $\Delta(x) = o(x^{1/2}\log^2x)$ has been deduced by Heath-Brown \cite{HB/1982} from the hypothesis $F(x,T)=o(T\log^2T)$ uniformly for $x^\eps\leq T \leq x$, and the sharper estimate $\Delta(x) = O(x^{1/2}\log^{5/4}x)$ is stated in Liu $\&$ Ye \cite{Li-Ye/2001} as a consequence of a quantitative version of the general pair-correlation conjecture by Rudnick $\&$ Sarnak \cite{Ru-Sa/1996}, specialized to the case of the Riemann zeta function. However, \cite{Li-Ye/2001} contains several inaccuracies (in part detected by Goldston \cite{Gol/2002} and Chan \cite{Cha/2003}), hence the results are not reliable although the general strategy is clear. From the point of view of the function $F(X,T,\tau)$, the approach in \cite{HB/1982} and \cite{Li-Ye/2001} may be formalized as follows.

\medskip
{\bf Theorem 3.} {\sl Let $x\to\infty$, $Z=x^{1/2}\log^2x$ and $U=U(x)\in[10,Z]$. Assume RH and suppose that $F(x,T,\tau)\ll T\log T$ for some $\tau=\tau(x)\in[1/Z,1]$, uniformly for $U\leq T \leq Z$. Then}
\[
\Delta(x) \ll  x^{1/2}(\log^2U + \tau^{1/2} \log^{3/2}x).
\]

\medskip
We omit the proof of Theorem 3 as it follows along the lines of Theorem 1 of Heath-Brown \cite{HB/1982} and Theorem 2 of Liu $\&$ Ye \cite{Li-Ye/2001}, using Lemma 6 below and partial summation. 

\medskip
{\bf Remarks.} {\bf 3.} The same argument allows a more general result, where the hypothetical bound on $F(x,T,\tau)$ is of type $O(T\log^a T)$ with $1\leq a \leq 2$, with suitable uniformity ranges. This produces a continuous chain of bounds for $\Delta(x)$ (exponents are linear in $a$), interpolating the classical result under RH and Theorem 3, in particular giving back Heath-Brown's bound.

{\bf 4.} Assume that the hypothesis in Theorem 3 holds with the choice
\begin{equation}
\label{1-13}
\tau = \min(1,\frac{\log^4U}{\log^3x}),
\end{equation}
so that
\begin{equation}
\label{1-14}
\Delta(x) \ll x^{1/2} \log^2U.
\end{equation}
Then, for example, the bound for $\Delta(x)$ in \cite{Li-Ye/2001} follows from \eqref{1-14} if the hypothesis in Theorem 3 holds with $U = \exp(\log^{5/8}x)$ and $\tau=1/\sqrt{\log x}$. Moreover, if such hypothesis holds with $U=\log\log x$ and $\tau=(\log\log\log x)^4/\log^3x$, then \eqref{1-14} gives
\[
\Delta(x) \ll x^{1/2} (\log\log\log x)^2,
\]
essentially Montgomery's conjectural upper bound for $\Delta(x)$ in \cite{Mon/1980}. Note that, apart from Heath-Brown \cite{HB/1982}, all these results require a very wide uniformity in $T$, entering the range $T\leq x^\eps$; however, the required value of $\tau$ is not very small. In Theorem 1 we have the opposite situation, i.e. $T$ lies in a plausible range but we require a much wider uniformity in $\tau$.

{\bf 5.} In view of the classical $\Omega$-results in prime number theory, \eqref{1-13} and \eqref{1-14} show that there are definite limitations to the uniformity ranges of the bound $F(x,T,\tau)\ll T\log T$. A similar remark could be applied to Theorem 1 as well (once such results are refined by replacing $X^\eps$ by the expected powers of $\log X$) in view of the oscillation results of Maier's \cite{Mai/1985} type for primes in short intervals. \fine

\medskip
We conclude observing that it would be desirable to have some numerical evidence towards bounds and uniformity ranges for $F(X,T,\tau)$, but such computations could be quite heavy.

\bigskip
{\bf Acknowledgements.} 
This research was partially supported by the grant PRIN2010-11 {\sl Arithmetic Algebraic Geometry and Number Theory}.
The first named author was partially supported by the Cariparo  
``Eccellenza'' grant \textsl{Differential Methods in Arithmetic, 
Geometry and Algebra}.

\bigskip
\section{Proof of Theorem 1}

\smallskip 
We assume from the beginning that $1 \leq h \leq X$ and $0\leq Y \leq X$. We start with the classical explicit formula, see Davenport \cite{Dav/2000}, in the form
\[
\psi (x) = x - \sum_{|\gamma|\le  X} \frac{x^\rho}{\rho} + O(\log^2X)
\] 
uniformly for $X\leq x \leq 2X$, thus obtaining
\begin{equation}
\label{2-1}
J(X,Y,h) \ll \int_X^{X+Y} \big|\sum_{|\gamma|\le  X} \frac{(x+h)^\rho - x^\rho}{\rho}\big|^2 \dx x+ YX^\eps.
\end{equation}
Let $1\leq U\le  X$ be a parameter to be chosen later on. We first evaluate the contribution of the terms in \eqref{2-1} with $\vert \gamma \vert \le  U$. Since $\rho=1/2+i\gamma$ we have
\[
\sum_{\vert \gamma \vert\le  U}  \frac{(x+h)^\rho-x^\rho}{\rho} = \int_{x}^{x+h} \sum_{\vert \gamma \vert\le  U}  u^{\rho-1}\ \dx u \ll x^{-1/2} \int_{x}^{x+h} \Bigl\vert \sum_{\vert \gamma \vert\le  U}  u^{i\gamma} \Bigr\vert\ \dx u,
\]
hence by the Cauchy-Schwarz inequality we have
\begin{equation}
\label{2-2}
\int_{X}^{X+Y}  \Big \vert \sum_{\vert \gamma \vert\le  U} \frac{(x+h)^\rho-x^\rho}{\rho}\Big \vert^2\dx x 
\ll \frac{h}{X} \int_{X}^{X+Y} \Big( \int_{x}^{x+h} \Bigl\vert \sum_{\vert \gamma \vert\le  U}  u^{i\gamma}
\Bigr\vert^2\dx u \Big) \dx x.
\end{equation}
We deal with the right hand side of  \eqref{2-2} first changing the order of integration and then adapting part of the arguments in Heath-Brown \cite{HB/1982}.  

\medskip
{\bf Lemma 1.} {\sl For $Y,h\geq 0$ and a continuous function $f(u)$ we have}
\[
\int_X^{X+Y}\Big(\int_x^{x+h} f(u)\dx u\Big)\dx x = \int_X^{X+h}\Big(\int_x^{x+Y} f(u)\dx u\Big)\dx x.
\]

\medskip
{\it Proof.} Let $x\geq X$ and $g(x) = \int_X^x f(u)\dx u$. Then the left hand side is
\[
\begin{split}
\int_X^{X+Y} \big(g(x+h)&-g(x)\big)\dx x = \int_{X+h}^{X+Y+h} g(x)\dx x - \int_{X}^{X+Y} g(x)\dx x \\
&= \int_{X+Y}^{X+Y+h} g(x)\dx x - \int_{X}^{X+h} g(x)\dx x = \int_X^{X+h} \big(g(x+Y)-g(x)\big)\dx x
\end{split}
\]
and the lemma follows. \fine

\medskip
We extend the definition of $\Sigma(X, T)$ in \eqref{1-2} writing for $v\in\R$
\[
\Sigma(X,T;v) =  \sum\limits_{|\gamma|  \le  T} X^{i\gamma} e^{i\gamma v}.
\]

\medskip
{\bf Lemma 2.} {\sl  For $V,T \geq1$ and $\tau\in(0,1]$ we have}
\[ 
\int_{\R} \vert \Sigma(V,T;v) \vert^2 e^{-2 \vert v \vert /\tau} \dx v = \tau F(V,T,\tau).
\]

\medskip
{\it Proof.} This is a variation on Lemma 3 of \cite{HB/1982}. Squaring out and using the substitution $v=y\tau$, the left hand side equals
\[
\tau \sum_{-T \le  \gamma , \gamma'\le  T} V^{i(\gamma - \gamma')} \int_{\R} e^{iy\tau(\gamma - \gamma')}
e^{-2 \vert y \vert } \dx y = \tau F(V,T,\tau),
\]
thanks to the Fourier transform formula
\[
\int_{\R} e^{iy\tau (\gamma - \gamma')} e^{-2 \vert y \vert } \dx y = \frac{4}{4+\tau^2(\gamma-\gamma')^2};
\]
see Lemma 3 of \cite{HB/1982}. \fine

\medskip
{\bf Lemma 3.} {\sl For $V,T \geq1$ and $\tau\in[0,1]$ we have}
\[
\int_{V}^{V(1+\tau)} \Bigl\vert  \sum_{|\gamma|\le  T} u^{i\gamma} \Bigr\vert^2  \dx u \ll V \tau F(V,T,\tau).
\]

\medskip
{\it Proof.} The case $\tau=0$ is trivial. By the substitution $u=Ve^{v}$, and observing that $e^{v} \le  e^3 e^{-2v/\tau}$ for $\tau\in (0,1]$ and $v\in [0,\log(1+\tau)]$, similarly as in (9) of \cite{HB/1982} we get that the left hand side equals
\[
\begin{split}
V \int_{0}^{\log(1+\tau)}  \vert \Sigma(V,T;v) \vert^2 e^{v}\dx v &\ll V \int_{0}^{\log(1+\tau)}  \vert \Sigma(V,T;v) \vert^2 e^{-2v/\tau}\dx v \\
&\ll  V\int_{\R}  \vert \Sigma(V,T;v) \vert^2 e^{- 2\vert v \vert/ \tau}\dx v =  V \tau F(V,T,\tau),
\end{split}
\]
thanks to Lemma 2. \fine

\medskip
Changing the order of integration in \eqref{2-2} by Lemma 1, and then applying Lemma 3 with $V=x, T=U$ and $\tau=Y/x$, from \eqref{2-2} we obtain
\begin{equation}
\label{2-3}
\int_{X}^{X+Y}  \Big \vert \sum_{\vert \gamma \vert\le  U} \frac{(x+h)^\rho-x^\rho}{\rho}\Big \vert^2\dx x \ll  \frac{h^2Y}{X} \max_{X\le  x \le  X+h} F(x, U, \frac{Y}{x}).
\end{equation}

\medskip 
Next we consider the contribution from the terms $U<\vert \gamma \vert \le  X$. We need a further lemma.

\medskip
{\bf Lemma 4.} {\sl For $V\geq 1$ and $\tau\in[0,1]$ we have}
\[
\int_{V}^{V(1+\tau)} \Bigl \vert \sum_{U< \vert \gamma \vert \le  X}  \frac{u^{i\gamma}}{\rho} \Bigr\vert^2 \dx u \ll 
V \tau \Big(\frac{F(V,X,\tau)}{X^2} + \frac{F(V,U,\tau)}{U^2} + \frac{1}{U^{1/2}}  \int_U^X F(V,t,\tau) \frac{\dx t}{t^{5/2}}\Big).
\]

\medskip
{\it Proof.} Analogously to p.93-94 of \cite{HB/1982}, applying partial summation to the inner sum, using the Cauchy-Schwarz inequality and applying Lemma 3 with $U\leq T\leq X$ we have
\begin{align*}
\int_{V}^{V(1+\tau)} & \Bigl \vert \sum_{U<\vert \gamma \vert\le X} \frac{u^{i\gamma}}{\rho} \Bigr \vert^2 \dx u 
\ll  \frac{1}{X^2} \int_{V}^{V(1+\tau)} \Bigl \vert \sum_{U<\vert \gamma \vert\le X} u^{i\gamma} \Bigr \vert^2 \dx u 
+ \int_{V}^{V(1+\tau)} \Bigl ( \int_U^X \Bigl \vert \sum_{U<\vert \gamma \vert\le  t} u^{i\gamma} \Bigr \vert  \frac{\dx t}{t^2}
\Bigr )^{2} \dx u \\
&\ll \frac{V \tau}{X^2} \Bigl(F(V,X,\tau) + F(V,U,\tau) \Bigr) + \int_{V}^{V(1+\tau)} \Big(\int_U^X \Big \vert \sum_{U<\vert \gamma \vert\le t}  u^{i\gamma} \Big \vert^2 \frac{\dx t}{t^{5/2}}\Big) \Big(\int_U^X \frac{\dx t}{t^{3/2}} \Big) \dx u \\
&\ll \frac{V \tau}{X^2} \Bigl( F(V,X,\tau) + F(V,U,\tau) \Bigr) + \frac{V \tau}{U^{1/2}}\int_U^X \Bigl(F(V,t,\tau) + F(V,U,\tau)\Bigr) \frac{\dx t}{t^{5/2}} \\
&\ll V \tau\Bigl(\frac{F(V,X,\tau)}{X^2} + \frac{F(V,U,\tau)}{U^2} +\frac{1}{U^{1/2}} \int_U^X F(V,t,\tau)\frac{\dx t}{t^{5/2}}\Bigr),
\end{align*}
which proves the lemma. \fine

\medskip
Writing
\[
\Bigl\vert\sum_{U<\vert \gamma \vert\le  X} \frac{(x+h)^\rho-x^\rho}{\rho}\Bigr\vert^2 \ll \Bigl \vert \sum_{U<\vert \gamma \vert\le  X} \frac{(x+h)^\rho}{\rho} \Bigr \vert ^2+ \Bigl \vert \sum_{U<\vert \gamma \vert\le X} \frac{x^\rho}{\rho} \Bigr \vert^2,
\]
applying Lemma 4 with $V=X+h$, $\tau=Y/(X+h)$ and with $V=X$, $\tau = Y/X$ we obtain
\begin{equation}
\label{2-4}
\begin{split}
&\int_{X}^{X+Y} \Big \vert \sum_{U<\vert \gamma \vert\le X}\frac{(x+h)^\rho-x^\rho}{\rho}\Big \vert^2\dx x \ll \\
& \hskip.5cm YX \Bigl\{ \frac{F(X+h,X,\frac{Y}{X+h})}{X^2} + \frac{F(X+h,U,\frac{Y}{X+h})}{U^2}  + \frac{1}{U^{1/2}}\int_{U}^{X} F\big(X+h,t,\frac{Y}{X+h}\big) \frac{\dx t}{t^{5/2}} \\
&\hskip1.5cm + \frac{F(X,X,\frac{Y}{X})}{X^2} + \frac{F(X,U,\frac{Y}{X})}{U^2} + \frac{1}{U^{1/2}}\int_{U}^{X} F\big(X,t,\frac{Y}{X}\big) \frac{\dx t}{t^{5/2}} \Bigr\}.
\end{split}
\end{equation}

\medskip
Assume now hypothesis $H(\eta)$ with $0<\eta<1$. Since we use the bound $F(X,T,\tau) \ll TX^\eps$ in the range $U\leq T \leq X$, from  \eqref{2-1}, \eqref{2-3} and \eqref{2-4} we may write 
\begin{equation}
\label{2-5}
\begin{split}
J(X,Y,h) &\ll  \frac{h^2Y}{X}F(X,U,\frac{Y}{X}) + YX \frac{F(X,U,\frac{Y}{X})}{U^2} + YX^\eps \\
&\ll  \frac{h^2YUX^\eps}{X} +\frac{YX^{1+\eps}}{U} + YX^\eps
\end{split}
\end{equation}
under the conditions imposed by hypothesis $H(\eta)$, namely
\begin{equation}
\label{2-6}
X^\eta \leq U \leq X \  \ \text{and} \ \ \frac{X^\eta}{U} \leq \frac{Y}{X} \leq 1.
\end{equation}
The optimal choice of $U$ in \eqref{2-5} is 
\begin{equation}
\label{2-7}
U=\frac{X}{h},
\end{equation}
thus getting the required bound for $J(X,Y,h)$ in \eqref{1-12}. Moreover, \eqref{2-6} and \eqref{2-7} give the uniformity conditions on $h$ and $Y$ in the first part of Theorem 1. When $\eta=0$ we require
\begin{equation}
\label{2-8}
X^\eps \leq U \leq X \  \ \text{and} \ \ 0 \leq \frac{Y}{X} \leq 1,
\end{equation}
therefore \eqref{2-7} and \eqref{2-8} give the uniformity conditions in the second part of Theorem 1.

\bigskip
\section{Comparing $F(X,T,\tau)$ and $\Sigma(X,T)$}

\smallskip
In this section we link $F(X,T,\tau)$ and $\Sigma(X,T)$, as required in Remark 1 of Section 1 and by the applications of $F(X,T,\tau)$ in Theorem 2.

\medskip
{\bf Lemma 5.} {\sl Let $\eps>0$, $X^\eps\leq T \leq X$ and $\tau T \leq  X^\eps$. Then}
\[ 
F(X,T,\tau) \ll T^{1+\eps}  + X^{\eps} \max_{t\leq T} \vert \Sigma(X,t) \vert^2.
\]

\medskip
{\it Proof.}  For $\tau=0$ the result is trivial, so we assume that $\tau>0$. Writing $\delta=(\tau/2)(1-\eps/2) \log T$ and recalling the definition of $\Sigma(X,T;v)$ after Lemma 1, thanks to the trivial bound $\Sigma(X,T;v)  \ll  T\log T$ we get
\[
\Bigl(\int_{-\infty}^{-\delta}+ \int_{\delta}^{+\infty} \Bigr) \vert \Sigma(X,T;v) \vert^2e^{-2 \vert v \vert /\tau}\dx v \ll   
\tau T^{1+\eps}.
\]
Hence in view of Lemma 2 with $V=X$ we obtain
\[ 
\tau F(X,T,\tau) = \int_{-\delta}^{\delta} \vert \Sigma(X,T;v) \vert^2 e^{-2 \vert v \vert /\tau} \dx v +O(\tau T^{1+\eps}).
\]
By partial summation we have
\[
\Sigma(X,T;v) =e^{ivT}\sum_{|\gamma|\le  T}X^{i\gamma}-iv \int_{0}^{T}\big( \sum_{\vert\gamma\vert\le  t} X^{i\gamma} \big) e^{ivt} \dx t\ll(1+\vert v \vert T)\max_{t\leq T}  \vert \Sigma(X,t) \vert,
\]
hence thanks to the hypothesis $\tau T \leq  X^\eps$ we get
\[
\begin{split}
\int_{-\delta}^{\delta} \vert \Sigma(X,T;v) \vert^2 e^{-2 \vert v \vert /\tau} \dx v &\ll \tau \log T (1+\tau^2 T^2 \log^2 T)
\max_{t\leq T}  \vert \Sigma(X,t) \vert^2 \\
&\ll \tau X^{\eps} 
\max_{t\leq T} \vert \Sigma(X,t) \vert^2.
\end{split}
\]
Lemma 5 is therefore proved. \fine

\medskip
In the opposite direction, since $F(X,T,0)=|\Sigma(X,T)|^2$ we have that hypothesis $H(0)$ implies Gonek's conjecture in the range $X^\eps\leq T \leq X$. We expect that a weaker hypothesis, namely the bound \eqref{1-10} uniformly for $X^\eps\leq T \leq X$ and $1/T\leq \tau\leq 1$, should already imply Gonek's conjecture in the range $X^\eps\leq T \leq X$. However, all we can prove at present is the following modified form of Lemma 1 of Heath-Brown $\&$ Goldston \cite{HB-Go/1984}. Writing
\begin{equation}
\label{3-1}
\Sigma(X,U,T)=\sum_{U<\vert \gamma \vert \leq T}  X^{i\gamma}
\end{equation}
and observing that $F(X,T,\tau)\geq 0$ thanks to Lemma 2, we have

\medskip
{\bf Lemma 6.} {\sl For $X,T\geq 2$, $0\leq U < T$ and $0\leq \tau\leq 1$ we have}
\[
\Sigma(X,U,T) \ll (1+T\tau)^{1/2}\, \max_{U\leq t\le  T} F(X,t,\tau)^{1/2}.
\]

\medskip
{\it Proof.} The case $\tau=0$ is trivial since $\vert \Sigma(X,0,T) \vert^2 = F(X,T,0)$, hence
\[
\vert \Sigma(X,U,T) \vert \leq \vert \Sigma(X,0,T) \vert + \vert \Sigma(X,0,U)\vert \ll \max_{U\leq t\le  T} F(X,t,0)^{1/2}.
\]
For $\tau>0$ we write
\[
\Sigma(X,U,T;v) =\sum_{U< \vert \gamma \vert \leq T}  X^{i\gamma} e^{iv\gamma}
\]
and apply the Sobolev-Gallagher inequality
\[
\vert f(0) \vert \leq \frac{1}{2\tau} \int_{-\tau}^{\tau}\vert f(v) \vert \dx v + \frac12  \int_{-\tau}^{\tau}\vert f'(v) \vert \dx v
\]
with $f(v)=\Sigma(X,U,T;v)^{2}$, thus obtaining
\begin{equation}
\label{3-2}
\vert \Sigma(X,U,T)\vert^{2} \ll \frac{1}{\tau} \int_{-\tau}^{\tau}\vert \Sigma(X,U,T;v) \vert^{2}\dx v +  \int_{-\tau}^{\tau}
\vert \Sigma(X,U,T;v) \vert \Bigl\vert \frac{\partial}{\partial v} \Sigma(X,U,T;v) \Bigr\vert  \dx v.
\end{equation}
By partial summation we have
\[
\Bigl\vert \frac{\partial}{\partial v} \Sigma(X,U,T;v) \Bigr\vert = \Bigl\vert \sum_{U<\vert \gamma \vert \leq T} \gamma X^{i\gamma} e^{iv\gamma} \Bigr\vert  \leq T \vert \Sigma(X,U,T;v) \vert + \Bigl\vert \int_{U}^{T} \Sigma(X,U,t;v) \dx t \Bigr\vert,
\]
hence by the Cauchy-Schwarz inequality we get
\begin{equation}
\label{3-3}
 \int_{-\tau}^{\tau} \Bigl\vert \frac{\partial}{\partial v} \Sigma(X,U,T;v) \Bigr\vert^{2} \dx v 
\ll T ^{2} \int_{-\tau}^{\tau}\vert\Sigma(X,U,T;v)\vert^{2}\dx v+T^2\max_{U\leq t\leq T} \int_{-\tau}^{\tau}\vert\Sigma(X,U,t;v)\vert^{2} \dx v.
\end{equation}
But, recalling that $\Sigma(X,0,T;v)=\Sigma(X,T;v)$, for $U\leq t \leq T$ Lemma 2 gives
\begin{equation}
\label{3-4}
\begin{split}
\int_{-\tau}^{\tau} \vert\Sigma(X,U,t;v)\vert^{2}\dx v&\ll \int_{-\infty}^{\infty} \left(\vert\Sigma(X,0,t;v)\vert^{2} +  \vert\Sigma(X,0,U;v)\vert^{2}\right) e^{-2\vert v \vert/\tau} \dx v \\
&\ll \tau \max_{U\leq t\leq T} F(X,t,\tau),
\end{split}
\end{equation}
therefore Lemma 6 follows by the Cauchy-Schwarz inequality from \eqref{3-2}, \eqref{3-3} and \eqref{3-4}. \fine

\bigskip
\section{Proof of Theorem 2} 

\smallskip
Let $\eps>0$ be arbitrarily small; in what follows the dependence on $\eps$ of the ranges of $h$ and $\eta$ will be written in a slightly weaker but cleaner form. We first show that the estimate for $J(X,Y,h)$ in Theorem 1, under $H(\eta)$ with $\eta>0$, implies
\begin{equation}
\label{4-1}
\psi(x+h)-\psi(x)=h+O(h^{2/3}x^{\eta/3+\eps}).
\end{equation}
Moreover, \eqref{4-1} is non-trivial (in the sense of Remark 2) for $x^{\eta + 5\eps} \leq h \leq x^{3/4 - \eta/2}$ and $ 0 < \eta < 1/2 - 4\eps$. The case $\eta=0$ will be treated by a direct application of $H(0)$, since we already observed in Remark 2 that both approaches give the same result in this case. We start noticing that if there exist $x_0\in[X,X+Y]$ and $\log X \leq K=o(h)$ such that $|\psi(x_0+h)-\psi(x_0)-h|\geq 5K$, then 
\begin{equation}
\label{4-2}
|\psi(x+h)-\psi(x)-h|\geq K \ \ \text{for every} \ \ x\in[X,X+Y]\cap \left[x_0-\frac{K}{\log X},x_0+\frac{K}{\log X}\right].
\end{equation}
Indeed, clearly both $\psi(x+h)$ and $\psi(x)$ may change at most by $\pm K$ with respect to their values at $x=x_0$, as $x$ runs over such an interval. Assuming that the above inequality holds for some $x_0\in[X,X+Y]$ with $K=h^a$ and $0<a<1$, from Theorem 1 with the choice $Y=hX^\eta$ we get
\begin{equation}
\label{4-3}
\frac{h^{3a}}{\log X} \ll J(X,Y,h) \ll h^2X^{\eta+\eps},
\end{equation}
since the cardinality of the $x$ in \eqref{4-2} is $\gg h^a/\log X$. Hence \eqref{4-1} follows, since \eqref{4-3} is contradictory if $h^a> ch^{2/3}X^{\eta/3+\eps}$ with a certain $c>0$, say. Actually, this argument gives \eqref{4-1} in a wider range of $h$, but a computation shows that the bound is non-trivial essentially only in the stated range.

\medskip
Next we show that a direct application of $H(\eta)$ gives
\begin{equation}
\label{4-4}
\psi(x+h)-\psi(x) - h \ll
\begin{cases}
h^{1/3} x^{1/6+\eta/3+\eps} & \text{if} \ 0<\eta<1 \\
h^{1/2} x^{\eps} & \text{if} \ \eta=0; 
\end{cases}
\end{equation}
\eqref{4-4} is non-trivial (again in the sense of Remark 2) for  $x^{1/4 + \eta/2 + 3 \eps} \leq h \leq x^{1 - \eta}$ if $0 < \eta < 1/2 - 2 \eps$ and for $x^{3 \eps} \leq h \leq x^{1 - \eps}$ if $\eta=0$. The proof of \eqref{4-4} is along the lines of Theorem 4 in Gonek \cite{Gon/1993}, especially in the case $\eta=0$. By the explicit formula we have
\begin{equation}
\label{4-5}
\psi(x+h) - \psi(x) = h - \sum_{|\gamma|\leq x/h} \int_x^{x+h}u^{\rho-1}\dx u  - \sum_{x/h < \vert \gamma \vert \leq x}\frac{(x+h)^\rho-x^{\rho}}{\rho} + O(\log^{2} x).
\end{equation}
Assume first that $0<\eta<1$ and recall \eqref{1-2}, \eqref{3-1} and the trivial bound $\Sigma(X,T)\ll T\log T$. Hence, applying Lemma 6 with $x\leq X\leq x+h$, $x^\eta\leq U\leq x/h$, $T=x/h$ and $\tau=x^\eta/U$, from hypothesis $H(\eta)$ we obtain
\begin{equation}
\label{4-6}
\begin{split}
\sum_{|\gamma|\leq x/h} \int_{x}^{x+h} u ^{\rho-1} \dx u &\ll hx^{-1/2}\big(\max_{x\leq u\leq x+h}|\Sigma(u,U)| + \max_{x\leq u\leq x+h}|\Sigma(u,U,\frac{x}{h})|\big) \\
&\ll hx^{-1/2+\eps}(U + \frac{x^{1+\eta/2}}{hU^{1/2}})  \ll h^{1/3} x^{1/6 + \eta/3 + \eps}
\end{split}
\end{equation}
with the optimal choice $U=(x^{2+\eta}/h^2)^{1/3}$, provided $x^\eta\leq x/h\leq x$. Note that the above conditions on $U$ are satisfied with such a choice. Let now $u=x$ or $u=x+h$. By partial summation and Lemma 6 with $X=u$, $U=x/h$, $x/h\leq T \leq x$ and $\tau =hx^{\eta-1}$, from hypothesis $H(\eta)$ we get
\begin{equation}
\label{4-7}
\sum_{x/h < |\gamma| \leq x} \frac{u^\rho}{\rho} \ll x^{-1/2} \vert \Sigma(u,\frac{x}{h},x) \vert
+ x^{1/2} \int_{x/h}^{x} \vert \Sigma(u,\frac{x}{h},t)\vert \frac{\dx t}{t^2} \ll h^{1/2} x^{\eta/2+\eps},
\end{equation}
again provided  $x^\eta\leq x/h\leq x$. From \eqref{4-6} and \eqref{4-7} we see that \eqref{4-5} becomes
\[
\psi(x+h) - \psi(x) - h \ll (h^{1/3} x^{1/6 + \eta/3} + h^{1/2} x^{\eta/2})x^\eps.
\]
A simple computation shows that the term $h^{1/3} x^{1/6 + \eta/3}$ dominates if $h\leq x^{1-\eta}$ and is non-trivial in the range stated after \eqref{4-4}. Moreover, the term $ h^{1/2} x^{\eta/2}$ dominates if $h\geq x^{1-\eta}$, but gives a trivial result in such a range. Thus the bound \eqref{4-4} follows if $\eta>0$. If $\eta=0$ the argument is similar but simpler, and reduces essentially to the proof of Theorem 4 of Gonek \cite{Gon/1993}. Indeed, the freedom in the choice of $\tau$ given by hypothesis $H(0)$ allows to choose $U=x^\eps$ and $\tau=0$ in the application of Lemma 6 to the first row of \eqref{4-6}, thus giving
\begin{equation}
\label{4-8}
 \sum_{|\gamma|\leq x/h} \int_x^{x+h}u^{\rho-1}\dx u \ll h^{1/2}x^\eps,
\end{equation}
provided $x^{3 \eps} \leq h \leq x^{1 - \eps}$. Hence the bound \eqref{4-4} follows from \eqref{4-8} and \eqref{4-7} in the case $\eta=0$, and is non-trivial in the stated range. Theorem 2 follows from \eqref{4-1} and \eqref{4-4}.


\vskip2cm

\vskip 1cm
\noindent
Alessandro Languasco, Dipartimento di Matematica, Universit\`a
di Padova, Via Trieste 63, 35121 Padova, Italy. \url{languasco@math.unipd.it}

\medskip
\noindent
Alberto Perelli, Dipartimento di Matematica, Universit\`a di Genova, via
Dodecaneso 35, 16146 Genova, Italy. \url{perelli@dima.unige.it}

\medskip
\noindent
Alessandro Zaccagnini, Dipartimento di Matematica e Informatica, Universit\`a di Parma, Parco
Area delle Scienze 53/a, 43124 Parma, Italy. \url{alessandro.zaccagnini@unipr.it}
\end{document}